\documentclass[12pt]{amsart}

\usepackage[psamsfonts]{amssymb}
 \newtheorem{thm}{Theorem}[section]
 \newtheorem{cor}[thm]{Corollary}
 \newtheorem{lem}[thm]{Lemma}
 \newtheorem{prop}[thm]{Proposition}
 \theoremstyle{definition}
 \newtheorem{defn}[thm]{Definition}
 \theoremstyle{remark}
 \newtheorem{rem}[thm]{Remark}
 \numberwithin{equation}{section}

\begin{document}

\newcommand{\D}{\mathbb{D}}
\newcommand{\C}{\mathbb{C}}
\newcommand{\N}{\mathbb{N}}
\newcommand{\R}{\mathbb{R}}
\newcommand{\dist}{\operatorname{dist}}

\renewcommand{\qedsymbol}{$\blacksquare$}

\title{Traces of H\"ormander algebras on discrete sequences}

\author[X. Massaneda]{Xavier Massaneda}

\address{%
Departament de Matem\`atica Aplicada i An\`alisi\\
Universitat  de Bar\-ce\-lo\-na, Gran Via 585\\
08071-Bar\-ce\-lo\-na\\
Spain}

\email{xavier.massaneda@ub.edu}

\thanks{Partially supported by the Picasso programme (Action Integr\'ee)
HF2006-0211. First and second authors supported by MEC grant
MTM2005-008984-C02-02 and CIRIT grant 2005-SGR 00611.}

\author[J. Ortega-Cerd\`a]{Joaquim Ortega-Cerd\`a}
\address{%
Departament de Matem\`atica Aplicada i An\`alisi\\
Universitat  de Bar\-ce\-lo\-na, Gran Via 585\\
08071-Bar\-ce\-lo\-na\\
Spain}
\email{jortega@ub.edu}

\author[M. Ouna\"{\i}es]{Myriam Ouna\"{\i}es}
\address{%
Institut de Recherche Math\'ematique Avanc\'ee\\
Universit\'e Louis Pasteur, 7 Rue Ren\'e Des\-car\-tes\\
67084 Strasbourg CEDEX\\
France}
\email{ounaies@math.u-strasbg.fr}


\date{\today}


\begin{abstract} 
We show that a discrete sequence $\Lambda$ of the complex plane is the union of
$n$ interpolating sequences for the H\"ormander algebras $A_p$ if and only if
the trace of  $A_p$ on $\Lambda$ coincides with the space of functions on
$\Lambda$ for which the  divided differences of order $n-1$ are uniformly
bounded. The analogous result holds in the unit disk for Korenblum-type
algebras. 
\end{abstract}

\maketitle

\section{Definitions and statement}

A function $p:\C\longrightarrow\R_+$,  is called a \emph{weight} if
\begin{itemize}
\item[(w1)] There is a constant $K>0$ such that $p(z)\ge K
\ln(1+|z|^2)$.
\item[(w2)]  There are  constants  $D_0>0$ and $E_0>0$ such that whenever $\vert
z-w\vert \le1 $ then 
\[p(z)\le D_0 p(w)+E_0.\] 
\end{itemize}

Let $H(\C)$ denote the space of all entire functions. We consider the algebra
\[
A_p=\Bigl\{f\in H(\C),\ \  \forall z\in \C,\  |f(z)| \le A\, e^{B p(z)}\hbox{
for some } A>0, B>0\Bigr\}. 
\]

Condition (w1) implies that $A_p$ contains the polynomials and (w2) that it is
closed under differentiation.

\begin{defn}
Given a discrete subset $\Lambda\subset\C$ we denote by  $A_p(\Lambda)$ the
space of sequences $\omega(\Lambda)=\{\omega(\lambda)\}_{\lambda\in \Lambda}$ of
complex numbers such that there are constants $A,B>0$ for which  
\[
 \vert \omega(\lambda) \vert\le A e^{Bp(\lambda)},\quad  \lambda\in \Lambda.
\]

We say that $\Lambda$ is an \emph{interpolating sequence for} $A_p$ when for
every sequence $\omega(\Lambda)\in A_p(\Lambda)$ there exists $f\in A_p$ such
that $f(\lambda)=\omega(\lambda)$, $\lambda\in\Lambda$. In terms of the
restriction operator 
\[ 
\begin{split}
\mathcal R_\Lambda: A_p & \longrightarrow A_p (\Lambda)\\
 f\, & \mapsto\; \{{f(\lambda)}\}_{\lambda\in \Lambda},
\end{split} 
\]
$\Lambda$ is  interpolating when  $\mathcal R_\Lambda(A_p)=A_p(\Lambda)$.
\end{defn}

\begin{defn}
Let $\Lambda$ be a discrete sequence in $\C$ and $\omega$ a function given on
$\Lambda$. The \emph{divided differences of  $\omega$} are defined by induction
as follows 
\[
\begin{split}
\Delta^0 \omega(\lambda_1)  &=\omega(\lambda_1)\ ,\\
\Delta^j\omega(\lambda_1,\ldots,\lambda_{j+1}) 
&=\displaystyle\frac{\Delta^{j-1}\omega(\lambda_2,\ldots,\lambda_{j+1})-\Delta^{
j-1}\omega(\lambda_1,\ldots,\lambda_j)}{\lambda_{j+1}-\lambda_1}\qquad j\geq
1.\\
\end{split}
\]

For any $n\in \N$, denote 
\[
\Lambda^n=\{(\lambda_1,\ldots,\lambda_n)\in
\Lambda\times\stackrel{\stackrel{n}{\smile}}{\cdots}\times \Lambda\; :\; 
\lambda_j\not=\lambda_k\ \textrm{if}\  j\not=k\},
\] 
and consider the set $X^{n-1}_p(\Lambda)$ consisting of the
functions in $\omega(\Lambda)$ with divided differences of order $n$ uniformly
bounded with respect to the weight $p$, i.e., such that for some $B>0$
\[
\sup_{(\lambda_1,\ldots,\lambda_n)\in \Lambda^n} \vert
\Delta^{n-1}\omega(\lambda_1,\ldots,\lambda_n)\vert
e^{-B[p(\lambda_1)+\cdots+p(\lambda_n)]}<+\infty\ .
\]
\end{defn}

\begin{rem}
It is clear that $X^n_p(\Lambda)\subset X^{n-1}_p(\Lambda)\subset\cdots\subset
X^0_p(\Lambda)=A_p(\Lambda)$.

To see this assume that $\omega(\Lambda)\in X^n_p(\Lambda)$, i.e., there exists
$B>0$ such that
\begin{multline*}
C:=\sup_{(\lambda_1,\dots,\lambda_{n+1})\in\Lambda^{n+1}}
 \left|\frac{\Delta^{n-1}\omega(\lambda_2,\dots,\lambda_{n+1})-\Delta^{n-1}
\omega(\lambda_1,\dots,\lambda_{n})}{\lambda_{n+1}-\lambda_1}\right|\times\\
\times e^{-B[p(\lambda_1)+\cdots+p(\lambda_{n+1})]}<\infty\ .
\end{multline*}
Then, given $(\lambda_1,\dots,\lambda_{n})\in\Lambda^{n}$ and taking
$\lambda_1^0,\dots,\lambda_n^0$ from a finite set (for instance the $n$ first
$\lambda^0_j\in\Lambda$ different of all $\lambda_j$) we have
\begin{multline*}
\Delta^{n-1}\omega(\lambda_1,\dots,\lambda_{n})=
\frac{\Delta^{n-1}\omega(\lambda_1,\dots,\lambda_{n})-\Delta^{n-1}
\omega(\lambda^0_1,\lambda_1,\dots,\lambda_{n-1})}{\lambda_{n}-\lambda^0_1}
(\lambda_{n}-\lambda_0^1)+\\
+\frac{\Delta^{n-1}
\omega(\lambda^0_1,\lambda_1,\dots,\lambda_{n-1})-\Delta^{n-1}
\omega(\lambda^0_2
,\lambda^0_1,\dots,\lambda_{n-2})}{\lambda_{n-1}-\lambda^0_2}(\lambda_{n-1}
-\lambda^0_2)+\cdots+\\
\frac{\Delta^{n-1}\omega(\lambda^0_{n-1},\dots,\lambda^0_1,\lambda_{1})-\Delta^{
n-1}\omega(\lambda^0_n,\dots,\lambda^0_1)}{\lambda_{1}-\lambda^0_n}(\lambda_{1}
-\lambda^0_n)+
\Delta^{n-1}\omega(\lambda^0_{n},\dots,\lambda^0_1)
\end{multline*}
Then a direct estimate and (w1) show that for some $B>0$ there is a constant
$K(\lambda^0_1,\dots,\lambda^0_n)$ such that
\begin{align*}
\left|\Delta^{n-1}\omega(\lambda_1,\dots,\lambda_{n})\right|&\leq C\left(
e^{B[p(\lambda^0_1)+\cdots+p(\lambda_n)]}+ \cdots 
+e^{B[p(\lambda^0_{n-1})+\cdots+p(\lambda_{1})]}\right)\\
&\leq K(\lambda^0_1,\dots,\lambda^0_n) e^{B[p(\lambda_1)+\cdots+p(\lambda_n)]},
\end{align*}
and the statement follows.

\end{rem}

 The main result of this note is modelled after Vasyunin's description of the
sequences $\Lambda$ in the unit disk such that the trace of the algebra of
bounded holomorphic functions $H^\infty$ on $\Lambda$ equals the space of
(hyperbolic) divided differences of order $n$ (see \cite{Vas83}, \cite{Vas84}).
The analogue in our context is the following.

\begin{thm}[Main Theorem]\label{main} The identity
$\mathcal R_\Lambda(A_p)= X^{n-1}_p(\Lambda)$ holds if and only if $\Lambda$ is
the union of $n$ interpolating sequences for $A_p$.
\end{thm}

For the most usual of these weights there exists a complete description of the
$A_p$-interpolating sequences, both in analytic and geometric terms. This is the
case for doubling and radial weights (see \cite[Corollary 4.8]{Br-Li}), or for
non-isotropic weights of the form $p(z)=|\operatorname{Im} z|+\log(1+|z|)$
(see \cite[Theorem 1]{MOO}).

With similar techiques it should be possible to extend this result to an
Hermite-type interpolation problem with multiplicites, along the lines of
\cite{Myriam}.

\section{General properties}

We begin by showing that one of the inclusions of Theorem~\ref{main} is
inmediate. 

\begin{prop}\label{incl}
For all $n\in\N$,  the inclusion $\mathcal R_\Lambda(A_p)\subset
X^{n-1}_p(\Lambda)$ holds.
\end{prop}

\begin{proof}
Let $f\in A_p$. Let us show  by induction on $j\ge1$ that, for certain
constants $A, B>0$
\[ 
\vert \Delta^{j-1} f(z_1,\ldots,z_j)\vert \le A 
e^{B[p(z_1)+\cdots+p(z_j)]}\qquad\textrm{for all $(z_1,\ldots,z_j)\in \C^j$}.
\]
As $f\in A_p$, we have $\vert \Delta^0 f(z_1)\vert =\vert f(z_1)\vert \le A
e^{Bp(z_1)}$.

Assume that the property is true for $j$ and let $(z_1,\ldots,z_{j+1})\in
\C^{j+1}$. Fix $z_1,\ldots,z_j$ and consider $z_{j+1}$ as the variable in the
function

\[
\Delta^j f(z_1,\ldots, z_{j+1}) =\frac{\Delta^{j-1}
f(z_2,\ldots,z_{j+1})-\Delta^{j-1}f(z_1,\ldots,z_j)}{z_{j+1}-z_1}.
\]
By the induction hypothesis, 
\begin{multline*}
\vert \Delta^{j-1} f(z_2,\ldots,z_{j+1})-\Delta^{j-1}f(z_1,\ldots,z_j)\vert \le
\\ 
\le A \bigl(e^{B[p(z_2)+\cdots+p(z_{j+1})]}+e^{B[p(z_1)+\cdots+p(z_j)]}\bigr)
\le 2Ae^{B[p(z_1)+\cdots+p(z_{j+1})]}.
\end{multline*}
Thus, if $\vert z_{j+1}-z_1\vert \ge 1$, we easily deduce the desired estimate.
For $\vert z_{j+1}-z_1\vert \le 1$, by the maximum principle and (w2):
\begin{align*}
\vert \Delta^j f(z_1,\ldots, z_{j+1})\vert  & \le 2A \sup_{\vert \xi-z_1\vert =
1} e^{B[p(z_1)+\cdots+p(z_j)+p(\xi)]}
\\ &  \le Ae^{(B+D_0)[p(z_1)+\cdots+p(z_j)+p(z_{j+1})]}.
\end{align*}
\end{proof}

\begin{defn}
A sequence $\Lambda$ is \emph{weakly separated} if there exist constants
$\varepsilon>0$ and $C>0$ such that the disks $D(\lambda, \varepsilon
e^{-Cp(\lambda)})$, $\lambda\in\Lambda$, are pairwise disjoint.
\end{defn}

\begin{rem}
If $\Lambda$ is weakly separated then $X^0_p(V)=X^n_p(V)$, for all $n\in\mathbb
N$. 

To see this it is enough to prove (by induction) that $X^0_p(\Lambda)\subset
X^n_p(\Lambda)$ for all $n\in\mathbb N$. For $n=0$ this is trivial. 
Assume now that $X^0_p(\Lambda)\subset X^{n-1}_p(\Lambda)$. Given
$\omega(\Lambda)\in X^0_p(\Lambda)$ we have
\begin{align*}
\left|\Delta^n\omega(\lambda_1,\dots,\lambda_{n+1})\right|&=\left|\frac{\Delta^{
n-1}(\lambda_2,\dots,\lambda_{n+1})-\Delta^{n-1}(\lambda_1,\dots,\lambda_n)}{
\lambda_{n+1}-\lambda_1}\right|\\
&\leq \frac{2A}{\varepsilon}e^{(B+C)[p(\lambda_1)+\cdots+p(\lambda_{n+1})]}\ .
\end{align*}

\end{rem}

\begin{lem}\label{equiv}
Let $n\geq 1$. The following assertions are equivalent:
\begin{itemize}
\item[(a)] $\Lambda$ is the union of $n$  weakly separated sequences,

\item[(b)] There exist constants $\varepsilon>0$ and $C>0$
such that
\[
 \sup_{\lambda\in\Lambda} \# [\Lambda\cap D(\lambda,\varepsilon
e^{-Cp(\lambda)})]\leq n\ .
\]

\item[(c)] $X^{n-1}_p(\Lambda)=X^n_p(\Lambda)$.
\end{itemize}

\end{lem}

\begin{proof}
(a) $\Rightarrow$(b). This is clear, by the weak separation.
 
(b) $\Rightarrow$(a). We proceed by induction on $j=1,\ldots,n$. For $j=1$, it
is again clear by the definition of weak separation. Assume  the property true
for $j-1$. Let $1\ge \varepsilon>0$ and $C>0$ be such that
$\sup_{\lambda\in\Lambda} \# [\Lambda\cap D(\lambda,\varepsilon
e^{-Cp(\lambda)})]\leq j$.
Put $\varepsilon'=e^{-E_0C}\varepsilon/2$ and $C'=D_0C$. By Zorn's Lemma, there
is a maximal subsequence  $\Lambda_1\subset\Lambda$ such that the disks
$D(\lambda, \varepsilon' e^{-C'p(\lambda)})$, $\lambda\in\Lambda_1$, are
pairwise disjoint. 
In particular $\Lambda_1$ is weakly separated.
For any  $\alpha \in \Lambda\setminus \Lambda_1$, there exists $\lambda\in
\Lambda_1$ such that
\[
D(\lambda, \varepsilon' e^{-C'p(\lambda)})\cap D(\alpha, \varepsilon'
e^{-C'p(\alpha)})\not=\emptyset,
\]
otherwise $\Lambda_1$ would not be maximal.
Then $\lambda\in D(\alpha,\varepsilon e^{-Cp(\alpha)})$, since
\[
\vert \lambda-\alpha\vert < \varepsilon' e^{-C' p(\lambda)}+\varepsilon'e^{-C'
p(\alpha)}< \varepsilon e^{-Cp(\alpha)},
\]
by (w2).
Thus $D(\alpha,\varepsilon e^{-Cp(\alpha)})$ contains at most $j-1$ points of
$\Lambda\setminus \Lambda_1$.
We use the induction hypothesis to conclude that $\Lambda\setminus \Lambda_1$ is
the union of $j-1$ weakly separated sequences and, by consequence, $\Lambda$ is
the union of $j$ weakly separated sequences. 

(b)$\Rightarrow$(c). It remains to see that $X^{n-1}_p(\Lambda)\subset
X^{n}_p(\Lambda)$.
Given $\omega(\Lambda)\in X^{n-1}_p(\Lambda)$ and points
$(\lambda_1,\ldots,\lambda_{n+1}) \in \Lambda^{n+1}$, we have to estimate
$\Delta^n \omega(\lambda_1,\ldots,\lambda_{n+1})$. Under the assumption (b), at
least one of these $n+1$ points  is not in the disk $D(\lambda_1, \varepsilon
e^{-Cp(\lambda_1)})$. Note that $\Lambda^n$ is invariant by permutation of the
$n+1$ points, thus we may assume that $\vert \lambda_1 -\lambda_{n+1}\vert \ge
\varepsilon e^{-Cp(\lambda_1)}$. Using the fact that $\omega(\Lambda)\in
X^{n-1}_p(\Lambda)$, there are constants $A,B>0$ such that
\[
\begin{split}
\vert \Delta^n \omega (\lambda_1,\ldots,\lambda_{n+1})\vert & \le \frac{\vert
\Delta^{n-1} \omega (\lambda_2,\ldots,\lambda_{n+1})\vert+\vert \Delta^{n-1}
\omega (\lambda_1,\ldots,\lambda_n)\vert}{\vert \lambda_1 -\lambda_{n+1}\vert}
\\
& \le A e^{B[p(\lambda_1)+\cdots+p(\lambda_{n+1})]}.
\end{split}
\] 

(c)$\Rightarrow$(b). We prove this by contraposition. Assume that for all
$C,\varepsilon>0$, there exists $\lambda\in \Lambda$ such that $\# [\Lambda\cap
D(\lambda,\varepsilon e^{-Cp(\lambda)})]>n$. Since $\Lambda$ has no
accumulation points, for any fixed $C>0$,  we can extract from
$\Lambda$ a weakly separated subsequence
$\mathcal{L}=\{\alpha^l\}_{l\in \N}$ such that $\#[(\Lambda\setminus
\mathcal{L})\cap D(\alpha^l,1/l\, e^{-Cp(\alpha^l)})]\geq n$ for all $l$. Let us
call $\lambda_1^l,\ldots,\lambda_n^l$ the points of $\Lambda\setminus
\mathcal{L}$ closest to $\alpha^l$,
arranged  by increasing distance.
In order to construct a sequence $\omega(\Lambda)\in X^{n-1}_p(\Lambda)\setminus
X^n_p(\Lambda)$, put 
\[
\begin{split}
\omega(\alpha^l) & = \prod_{j=1}^{n-1} (\alpha^l-\lambda_j^l),\,  \hbox{for all
} \alpha^l\in \mathcal{L}\\
\omega(\lambda) & = 0 \hbox{ if } \lambda\in \Lambda\setminus \mathcal{L}.
\end{split}
\]
 
To see that $\omega(\Lambda)\in X^{n-1}(\Lambda)$ let us estimate 
$\Delta^{n-1}\omega(\lambda_1,\ldots,\lambda_n)$ for any given  vector
$(\lambda_1,\ldots,\lambda_n)\in \Lambda^n$. We don't need to consider the case
where the points are distant, thus, as $\mathcal{L}$ is weakly separated, we may
assume that at most one of the points is in $\mathcal{L}$. On the other hand, it
is clear that 
$\Delta^{n-1}\omega(\lambda_1,\ldots,\lambda_n)=0$ if all the points are in
$\Lambda\setminus \mathcal{L}$. Then, taking into account that $\Delta^{n-1}$ is
invariant by permutation, we will only consider the case where  $\lambda_n$ is
some $\alpha^l \in \mathcal{L}$ and  $\lambda_1,\ldots,\lambda_{n-1}$ are in
$\Lambda\setminus \mathcal{L}$. In that case, 
\[ 
\vert \Delta^{n-1}\omega(\lambda_1,\ldots,\lambda_{n-1},\alpha^l)\vert =\vert
\omega(\alpha^l)\vert \prod_{k=1}^{n-1} \vert \alpha^l-\lambda_k^l\vert^{-1}
 \le 1,
 \]
 as desired.

On the other hand, a similar computation yields
\[ 
\begin{split}
\vert \Delta^n \omega(\lambda_1^l,\ldots,\lambda_n^l,\alpha^l)\vert   =\vert
\omega(\alpha^l)\vert \prod_{k=1}^n\vert\alpha^l-\lambda_k^l\vert^{-1}
 =\vert \alpha^l-\lambda_n^l\vert^{-1} \ge l e^{Cp(\alpha^l)}.
\end{split}
\]
Using (w2), for any constant $B>0$, and choosing $C=B(nD_0+1)$, we have
\[
\vert \Delta^n \omega(\lambda_1^l,\ldots,\lambda_n^l,\alpha^l)\vert
e^{-B(p(\lambda_1^l)+\cdots+p(\lambda_n^l)+p(\alpha^l))}\ge l
e^{-BnE_0}\rightarrow +\infty.
\]
We finally conclude that $\omega(\Lambda) \notin X^n_p(\Lambda)$.
\end{proof}

\begin{cor}\label{weak}
If $\Lambda$ is an interpolating sequence, then it is weakly separated.
\end{cor}

\begin{proof}
If $\Lambda$ is an interpolating sequence, then $\mathcal R_\Lambda(A_p)=
X^0_p(\Lambda)$. On the other hand, by Proposition \ref{incl}, $\mathcal
R_\Lambda(A_p)\subset X^1_p(\Lambda)$. Thus  $X^0_p(\Lambda)= X^1_p(\Lambda)$.
We conclude by the preceding lemma applied to the particular case $n=1$.
\end{proof}

\begin{lem}\label{cover}
Let $\Lambda_1,\ldots,\Lambda_n$ be weakly separated sequences. 
There exist  positive constants $a,b$, $B_1,B_2$ and $\varepsilon>0$, a
subsequence $\mathcal{L}\subset \Lambda_1\cup\cdots\cup \Lambda_n$ and disks
$D_\lambda=D(\lambda, r_\lambda)$, $\lambda\in \mathcal{L}$, such that
\begin{itemize}
\item[(i)]  $\Lambda_1\cup\cdots\cup \Lambda_n\subset  \cup_{\lambda\in
\mathcal{L}} D_\lambda$ 
\item[(ii)]   $a\varepsilon e^{-B_1p(\lambda)}\le r_\lambda\le b \varepsilon
e^{-B_2p(\lambda)}$ for all $\lambda \in \mathcal{L}$
\item[(iii)]  $\dist(D_\lambda,D_{\lambda'})
\ge a\varepsilon e^{-B_1p(\lambda)}$ for all $\lambda,\lambda' \in \mathcal{L}$,
$\lambda\neq \lambda'$.
\item[(iv)] $\#(\Lambda_j\cap D_\lambda)\leq 1$ for all $j=1,\dots,n$ and
$\lambda\in \mathcal{L}$. 
\end{itemize}
\end{lem}

\begin{proof}
Let $0<\varepsilon<1$ and $C>0$ be constants such that 
\begin{equation}\label{separation}
\vert \lambda-\lambda'\vert \ge \varepsilon e^{-C/D_0(p(\lambda)-E_0)}, \qquad
\forall\lambda,\lambda'\in \Lambda_j,\ \lambda\neq\lambda', \qquad \forall
j=1,\ldots, n\ ,
\end{equation}
where $D_0\ge 1$ and $E_0\ge 0$ are given by (w2).

We will proceed by induction on $k=1,\ldots,n$ to show the existence of a
subsequence $\mathcal{L}_k\subset \Lambda_1\cup\cdots\cup \Lambda_{k}$ and
constants $C_k\ge C$, $B_k\ge 0$ such that:
\begin{align*}
&(i)_k\quad \Lambda_1\cup\cdots\cup \Lambda_k \subset  \cup_{\lambda\in
\mathcal{L}_k}
D(\lambda, R_\lambda^k), \\
&(ii)_k\quad 2^{-3k}e^{-C_kp(\lambda)-B_k}\varepsilon \le R_\lambda^k\le
\varepsilon e^{-Cp(\lambda)}\sum_{j=0}^{k-1}2^{-(3j+2)}\le
2/7e^{-Cp(\lambda)}\varepsilon, \\
&(iii)_k\quad \dist(D(\lambda, R_\lambda^k),D(\lambda', R_{\lambda'}^k))
\ge 2^{-3k} \varepsilon e^{-C_kp(\lambda)-B_k}\, \textrm{for any 
$\lambda,\lambda'\in \mathcal{L}_k$, $\lambda\neq\lambda'$}.
\end{align*}

The constants $C_k$ and $B_k$ are chosen, in view of (w2), so that $C_k
p(\lambda)+B_k\le C_{k+1} p(\lambda')+B_{k+1}$ whenever $\vert
\lambda-\lambda'\vert \le 1$.

Then it suffices to chose $\mathcal{L}=\mathcal{L}_n$, $r_\lambda=R_\lambda^n$,
$a=e^{-B_n}2^{-3n}$, $b=2/7$, $B_1=C_n$ and $B_2=C$. 
As $r_\lambda< e^{-Cp(\lambda)}\varepsilon $, it is clear that
$D(\lambda,r_\lambda)$ contains at most one point of each  $\Lambda_j$, hence
the lemma follows.

For $k=1$, the property is clearly verified with  $\mathcal{L}_1=\Lambda_1$ and
$R_\lambda^1=e^{-Cp(\lambda)}\varepsilon /4 $. 

Assume the property true for $k$ and split $\mathcal{L}_k=\mathcal{M}_1\cup
\mathcal{M}_2$ and $\Lambda_{k+1}=\mathcal{N}_1\cup \mathcal{N}_2$, where 
\[
\begin{split}
\mathcal{M}_1 & =\{\lambda\in \mathcal{L}_k\ :\ D(\lambda,
R_\lambda^k+2^{-3k-2}\varepsilon e^{-C_k p(\lambda)-B_k} )\cap
\Lambda_{k+1}\not= \emptyset\}, \\
\mathcal{N}_1 & =\Lambda_{k+1}\cap\bigcup_{\lambda\in \mathcal{L}_k} 
D(\lambda, R_\lambda^k+2^{-3k-2}\varepsilon e^{-C_kp(\lambda)-B_k}), \\
\mathcal{M}_2 & = \mathcal{L}_k\setminus \mathcal{M}_1,\\
\mathcal{N}_2 & = \Lambda_{k+1}\setminus \mathcal{N}_1.
\end{split}
\]
Now, we put $\mathcal{L}_{k+1}=\mathcal{L}_k\cup \mathcal{N}_2$ and define the
radii $R_\lambda^{k+1}$ as follows: 
\[
R_\lambda^{k+1}  =
\begin{cases}
R_\lambda^k+2^{-3k-2} \varepsilon e^{-C_kp(\lambda)-B_k}\ &\textrm{if}\
\lambda\in  \mathcal{M}_1,   \\
R_\lambda^k \ &\textrm{if}\  \lambda\in  \mathcal{M}_2, \\
 2^{-3k-3} \varepsilon e^{-C_{k+1}p(\lambda)-B_{k+1}} \ &\textrm{if}\ 
\lambda\in  \mathcal{N}_2.
\end{cases}
\]
It is clear that 
 \[
\Lambda_1\cup\cdots\cup \Lambda_{k+1}\subset  \bigcup_{ \lambda\in
\mathcal{L}_{k+1}}
D(\lambda, R_\lambda^{k+1})
\]
and, by the induction hypothesis,
\[
2^{-3k-3}\varepsilon  e^{-C_{k+1}p(\lambda)+B_{k+1})} \le  R_\lambda^{k+1}\le
\varepsilon e^{-Cp(\lambda)}\sum_{j=0}^k 2^{-3j-2}\le 2/7\varepsilon
e^{-Cp(\lambda)}.
\]

In order to prove $(iii)_k$ take now  $\lambda,\lambda'\in\mathcal{L}_{k+1}$,
$\lambda\neq\lambda'$. We will verify that 
\[
\dist(D(\lambda,R_\lambda^{k+1}), D(\lambda',R_{\lambda'}^{k+1}))=
\vert \lambda-\lambda' \vert-R_\lambda^{k+1}-R_{\lambda'}^{k+1} \ge 2^{-3k-3} 
\varepsilon
e^{-C_{k+1}p(\lambda)-B_{k+1}} \]
by considering different cases.

If  $\lambda,\lambda'\in\mathcal{L}_k$ and $p(\lambda)\le p(\lambda')$, then
\begin{align*}
\dist(D(\lambda,R_\lambda^{k+1}), D(\lambda',R_{\lambda'}^{k+1}))&
\ge \vert \lambda-\lambda'\vert-R_\lambda^k-R_{\lambda'}^k- 2^{-3k-1}
\varepsilon e^{-C_kp(\lambda)-B_k} \\
&\ge 2^{-3k-1}\varepsilon e^{-C_kp(\lambda)-B_k} .
\end{align*}

Assume now $\lambda,\lambda'\in\mathcal{N}_2$ and $p(\lambda)\le p(\lambda')$.
Condition \eqref{separation} implies $\vert \lambda-\lambda'\vert\ge \varepsilon
e^{-Cp(\lambda)}$, hence
 \[
\dist(D(\lambda,R_\lambda^{k+1}), D(\lambda',R_{\lambda'}^{k+1}))
\ge (1-2^{-3k-2})  \varepsilon  e^{-Cp(\lambda)}.
\]

If $\lambda \in \mathcal{M}_1$ and ${\lambda'}\in \mathcal{N}_2$ there exists
$\beta\in \mathcal{N}_1 $ such that $\vert \lambda-\beta\vert\le
R_\lambda^{k+1}$. There is no restriction in assuming that $\vert
\lambda-\lambda'\vert\le 1$. Then, using \eqref{separation} on
$\beta,\lambda'\in\Lambda_{k+1}$,  we have
\[
\vert \lambda-\lambda'\vert \ge \vert \beta-\lambda' \vert-\vert
\lambda-\beta\vert\ge \varepsilon e^{-C/D_0(p(\beta)-E_0)}-R_\lambda^{k+1}\ge
\varepsilon e^{-Cp(\lambda)}-R_\lambda^{k+1}\ .
\]
The definition of $R_{\lambda'}^{k+1}$ together with the estimate
$R_{\lambda}^{k+1}\leq 2/7 \varepsilon e^{-Cp(\lambda)}$ yield 
\begin{multline*}
\dist(D(\lambda,R_\lambda^{k+1}), D(\lambda',R_{\lambda'}^{k+1})) 
\ge\varepsilon e^{-Cp(\lambda)} -2 R_{\lambda}^{k+1}-R_{\lambda'}^{k+1}\\
\ge \varepsilon e^{-Cp(\lambda)}-2R_\lambda^k-2^{-3k-1} \varepsilon 
e^{-C_{k}p(\lambda)-B_{k}}-2^{-3k-3} \varepsilon 
e^{-C_{k+1}p(\lambda')-B_{k+1}}\\
\ge    \varepsilon e^{-Cp(\lambda)}-\frac 47 \varepsilon e^{-Cp(\lambda)}
-2^{-3k} \varepsilon e^{-C_{k}p(\lambda)-B_{k}}
\ge \varepsilon e^{-Cp(\lambda)}(3/4-2^{-3k}) ,
\end{multline*}
as required.

Finally, if $\lambda \in \mathcal{M}_2$ and ${\lambda'}\in \mathcal{N}_2$,
again, assuming that $\vert \lambda-\lambda'\vert\le 1$, we have
\begin{align*}
\dist(D(\lambda,R_\lambda^{k+1}), D(\lambda',R_{\lambda'}^{k+1}))  &=\vert
\lambda-\lambda'\vert-R_\lambda^k-
2^{-3k-3} \varepsilon e^{-C_{k+1}p(\lambda')-B_{k+1}} \\
 &\ge 2^{-3k-2} \varepsilon e^{-C_kp(\lambda)-B_k} -2^{-3k-3} \varepsilon
e^{-C_{k}p(\lambda)-B_{k}}\\
 &\ge 2^{-3k-3} \varepsilon e^{-C p(\lambda)}   .
\end{align*}
\end{proof}

\section{Proof of Theorem~\ref{main}. Necessity}\label{necessity}

Assume $\mathcal R_\Lambda(A_p)= X^{n-1}_p(\Lambda)$, $n\ge 2$. Using
Proposition \ref{incl}, we have $X^{n-1}_p(V)=X^n_p(V)$, and by Lemma
\ref{equiv} we deduce that 
$\Lambda=\Lambda_1\cup \cdots\cup \Lambda_n$, where $\Lambda_1,\ldots,\Lambda_n$
are weakly separated sequences. We want to show that each $\Lambda_j$ is an
interpolating sequence. 

Let $\omega(\Lambda_j)\in A_p(\Lambda_j)= X^0_p(\Lambda_j)$.  
Let $\cup_{\lambda\in \mathcal{L}}D_\lambda$ be the covering of $\Lambda$ given
by Lemma \ref{cover}. 
We extend $\omega(\Lambda_j)$ to a sequence $\omega(\Lambda)$ which is constant
on each $D_\lambda\cap\Lambda_j$ in the following way:
\[
\omega_{|D_\lambda\cap\Lambda }=
\begin{cases}
0\quad &\textrm{if $D_\lambda\cap\Lambda_j=\emptyset$}\\
\omega(\alpha)\quad &\textrm{if $D_\lambda\cap\Lambda_j=\{\alpha\}$}\ .
\end{cases}
\]
We verify by induction that the extended sequence is in $X^{k-1}_p(\Lambda)$ for
all $k$.  It is clear that it belongs to $X^0_p(\Lambda)$. Assume that
$\omega\in X^{k-2}_p(\Lambda)$ and consider $(\alpha_1,\ldots,\alpha_k)\in
\Lambda^k$. If all the points are in the same $D_\lambda$ then
$\Delta^{k-1}\omega(\alpha_1,\ldots,\alpha_k)=0$, so we may assume that
$\alpha_1\in D_{\lambda}$ and $\alpha_k\in D_{\lambda'}$ with $\lambda
\not=\lambda'$. Then  we have
 \[
 \vert \alpha_1-\alpha_k \vert \ge a\varepsilon e^{-B_1p(\lambda)},
 \]
by Lemma~\ref{cover} (iii). With this and the induction hypothesis it is clear
that for certain constants $A,B>0$ 
\begin{align*}
\vert \Delta^{k-1}\omega(\alpha_1,\ldots,\alpha_k)\vert&=\left|\frac{
\Delta^{k-2}\omega(\alpha_2,\ldots,\alpha_k)-
\Delta^{k-2}\omega(\alpha_1,\ldots,\alpha_k)}{\alpha_1-\alpha_k}\right|\\
&\le A  e^{B[p(\alpha_1)+\cdots+p(\alpha_k)]}.
\end{align*}
In particular $\omega(\Lambda)\in X^{n-1}_p(\Lambda)$, and by assumption, there
exist $f\in A_p$ interpolating  the values $\omega(\Lambda)$. In particular $f$
interpolates  $\omega(\Lambda_j)$.
 
\section{Proof of Theorem~\ref{main}. Sufficiency}\label{sufficiency}

According to Proposition~\ref{incl} we only need to see that
$X^{n-1}_p(\Lambda)\subset \mathcal R_\Lambda(A_p)$. 
 
Before going further, let us recall the following facts about interpolation in
the spaces $A_p$.

\begin{lem}\cite[Lemma 2.2.6]{Br-Ga}\label{unif}
Let $\Gamma$ be an $A_p$-interpolating sequence. Then:
\begin{itemize}
\item[(i)] For all $A,B>0$, there exist constants $A',B'>0$ such that for all
sequences $\omega\in A_p(\Gamma)$ with 
$\displaystyle \sup_{\gamma\in \Gamma}\vert \omega(\gamma)\vert
e^{-Bp(\gamma)}\le A$ there exists $f\in A_p$ with $\displaystyle \sup_{z}\vert
f(z)\vert e^{-B'p(z)}\le A'$ and $f(\gamma)=\omega(\gamma)$ for all $\gamma\in
\Gamma$..

\item[(ii)] There exists a constant $C>0$ such that $\sum_{\gamma \in \Gamma}
e^{-Cp(\gamma)}<\infty$.

\end{itemize}
\end{lem}

Applying (i) to the sequences
$\omega_\gamma=\{\delta_{\gamma,\gamma'}\}_{\gamma'\in\Gamma}$ it is easy to
deduce that $\Gamma$ is weakly separated. Property (ii) is just a consequence of
the weak separation and properties (w1) and (w2).

Assume thus that $\Lambda=\Lambda_1\cup \cdots\cup \Lambda_n$ where $\Lambda_1,
\dots,\Lambda_n$ are interpolating sequences. Recall that each $\Lambda_j$ is
weakly separated (Corollary~\ref{weak}). Consider also the covering of $\Lambda$
given by Lemma~\ref{cover}.

\begin{lem}\label{pick}
There exist constants $A,B>0$ and a sequence  $\{F_\lambda\}_{\lambda\in
\mathcal{L}}\subset A_p$ such that:
\begin{align*}
&F_\lambda(\alpha)=
\begin{cases}
1\ &\textrm{ if }\ \alpha\in  \Lambda \cap D_\lambda\\
0  &\textrm{ if }\ \alpha\in  \Lambda\cap D_{\lambda'},\ \lambda'\not=\lambda 
\end{cases}\\
& \vert F_{\lambda}(z)\vert \le A e^{B(p(\lambda)+p(z))}\quad\textrm{for all
$z\in \C$}.\\
\end{align*}
\end{lem}

\begin{proof}
Fix $\lambda\in \mathcal{L}$ and define $\omega(\Lambda)$ by
\[
\omega(\alpha)=
\begin{cases}
\prod\limits_{\beta\in \Lambda \cap D_\lambda}(\alpha-\beta)^{-1}\quad 
&\textrm{ if }\ \alpha\notin \Lambda\cap D_\lambda  \\ 
\ 0 &\textrm{ if }\ \alpha\in \Lambda\cap D_\lambda.
\end{cases}
\]
By Lemma \ref{cover} (iii), we have $\vert \alpha-\beta\vert \ge c\varepsilon
e^{-Cp(\alpha)}$ whenever $\alpha\notin \Lambda \cap D_\lambda$, $\beta\in
\Lambda \cap D_\lambda$. 
Since $\#( \Lambda\cap D_\lambda)\leq n$ we deduce that
\[
\vert \omega(\alpha)\vert \le (c\varepsilon)^{-n}e^{nCp(\alpha)}
\]
Recall that $\Lambda_j$ is an interpolating sequence for all  $j=1,\ldots,n$, 
thus there exist a $n$-indexed sequence
$\{f_{\lambda,j}\}_{\lambda\in \mathcal{L},j
\in [[1,n]]}\subset A_p$ such that for all $z\in \C$,
\begin{align*}
\vert f_{\lambda,j}(z)\vert &\le A e^{Bp(z)} \\
f_{\lambda,j}(\alpha)&=\prod_{\beta\in \Lambda \cap
D_\lambda}(\alpha-\beta)^{-1}  \hbox{ if } \alpha\notin  \Lambda_j \cap
D_\lambda,
\end{align*}
with the constants $A$ and $B$ independent of $\lambda$  (see
Lemma~\ref{unif}(i) ).

The sequence of functions $\{F_\lambda\}_{\lambda\in \mathcal{L}}$ defined by  
\[ 
F_\lambda(z)=\prod_{j=1}^n\left[ 1-\prod_{\beta\in \Lambda \cap
D_\lambda}(z-\beta)f_{\lambda,j}(z)\right]
\]
has the desired properties.
\end{proof}

\begin{lem}\label{pick2}
For all $D>0$, there exist  $D'>0$ and a sequence  $\{G_\lambda\}_{\lambda\in
\mathcal{L}}\subset A_p$ such that:
\begin{align*}
&G_\lambda(\alpha)=
e^{Dp(\lambda)} \ \textrm{ if }\ \alpha\in  \Lambda \cap D_\lambda.\\
& \vert G_{\lambda}(z)\vert \le A e^{Bp(\lambda)}e^{D'p(z)}\quad\textrm{for all
$z\in \C$},\\
\end{align*}
where $A,B>0$ do not depend on $D$. 
\end{lem}

\begin{proof}
In this proof $D'$ denotes a constant depending on $D$ but not on $\lambda$, and
its actual value may change from one occurrence to the other.

Let $\lambda\in \mathcal{L}$. Assume, without loss of generality, that
$D_\lambda\cap \Lambda_j=\{\alpha_{\lambda,j}\}$
for all $j$.
As $\Lambda_1$ is an interpolating sequence and $e^{Dp(\lambda)}\le Ae^{D'
p(\alpha_{\lambda,1})}$, by Lemma \ref{unif}(i) there exists a sequence
$\{h_{\lambda,1}\}_\lambda \subset A_p$ such that 
\begin{align*}
&h_{\lambda,1}(\alpha_{\lambda,1})=e^{Dp(\lambda)}\\
& \vert h_{\lambda,1}(z)\vert \le A e^{D'p(z)}\quad\textrm{for all $z\in \C$}.\\
\end{align*}
Setting $H_{\lambda,1}(z)=h_{\lambda,1}(z)$, we have
$H_{\lambda,1}(\alpha_{\lambda,1})=e^{Dp(\lambda)}$.
Now, as $\Lambda_2$ is $A_p$-interpolating and 
\[ 
\frac{\vert e^{Dp(\lambda)}-H_{\lambda,1}(\alpha_{\lambda,2})\vert }{\vert
\alpha_{\lambda,2}-\alpha_{\lambda,1}\vert }=\frac{\vert
H_{\lambda,1}(\alpha_{\lambda,1})-H_{\lambda,1}(\alpha_{\lambda,2})\vert }{\vert
\alpha_{\lambda,2}-\alpha_{\lambda,1}\vert }\le Ae^{D' p(\alpha_{\lambda,2})},
\]
 there exists a sequence $\{h_{\lambda,2}\}_\lambda \subset A_p$ such that 
\begin{align*}
&h_{\lambda,2}(\alpha_{\lambda,2})=\frac{e^{Dp(\lambda)}-H_{\lambda,1}(\alpha_{
\lambda,2})}{\alpha_{\lambda,2}-\alpha_{\lambda,1}}\\
& \vert h_{\lambda,2}(z)\vert \le A e^{D'p(z)}\quad\textrm{for all $z\in \C$}.\\
\end{align*}
Setting
$H_{\lambda,2}(z)=h_{\lambda,1}(z)+h_{\lambda,2}(z)(z-\alpha_{\lambda,1})$. We
have 
\[ 
H_{\lambda,2}(\alpha_{\lambda,1})=H_{\lambda,2}(\alpha_{\lambda,2})=e^{
Dp(\lambda)}.
\]
We proceed by induction to construct a sequence of functions
$\{h_{\lambda,k}\}_\lambda \subset A_p$
such that
\begin{align*}
&h_{\lambda,k}(\alpha_{\lambda,k})=\frac{e^{Dp(\lambda)}-H_{\lambda,k-1}(\alpha_
{\lambda,k})}{(\alpha_{\lambda,k}-\alpha_{\lambda,1})\cdots(\alpha_{\lambda,k}
-\alpha_{\lambda,k-1})}\\
& \vert h_{\lambda,k}(z)\vert \le A e^{D'p(z)}\quad\textrm{for all $z\in \C$}.
\end{align*}

Then the function defined by
$H_{\lambda,k}(z)=H_{\lambda,k-1}(z)+h_{\lambda,k}(z)(z-\alpha_{\lambda,1}
)\cdots (z-\alpha_{\lambda,k-1})$ verifies 
\[ 
H_{\lambda,k}(\alpha_{\lambda,1})=\cdots=
H_{\lambda,k}(\alpha_{\lambda,k})=e^{Dp(\lambda)}.
\]
Finally, we set $G_\lambda=H_{\lambda,n}$. 
\end{proof}

To proceed with the proof of the inclusion $X^{n-1}_p(\Lambda)\subset
\mathcal{R}_\Lambda(A_p)$, let $\omega(\Lambda)\in X^{n-1}_p(\Lambda)$.

Fix $\lambda\in \mathcal{L}$ and let $\Lambda\cap
D_\lambda=\{\alpha_1,\ldots,\alpha_k\}$, $k\le n$.
We first consider a polynomial interpolating the values
$\omega(\alpha_1),\ldots, \omega(\alpha_k)$:
\[ 
P_\lambda(z)=\Delta^0\omega(\alpha_1)+\Delta^1\omega(\alpha_1,
\alpha_2)(z-\alpha_1)+\cdots+\Delta^{k-1}\omega(\alpha_1,\ldots,\alpha_k)\prod_{
j=1}^{k-1}(z-\alpha_j).
\]
Notice that $P_\lambda\in A_p$, since $\omega(\Lambda)\in X^{n-1}_p(\Lambda)$
and by properties (w1) and (w2) we have 
\[
|P_\lambda(z)|\le A\vert z\vert^k e^{B[p(\alpha_1)+\cdots+p(\alpha_k)]}\le A
e^{B'[p(z)+p(\lambda)]}\ .
\]

Now, define
\[
f=\sum_{\lambda\in \mathcal{L}} F_\lambda G_\lambda P_\lambda e^{-Dp(\lambda)},
\]
where $D$ is a large constant to be chosen later on.

By the preceding estimates on $G_\lambda$ and $P_\lambda$, there exist 
constants $A,B>0$ not depending on $D$ and a constant $D''>0$ such that, for all
$z\in \C$, we have 
\[
\vert f(z)\vert \le A e^{D''p(z)}\sum_{\lambda\in \mathcal{L}}
e^{(B-D)p(\lambda)}.
\]
In view of Lemma~\ref{unif} (ii), taking $D=B+C$, the latter sum converges and
$f\in A_p$.

To verify that $f$ interpolates $\omega(\Lambda)$, let $\alpha\in\Lambda$ and
let $\lambda$ be the (unique) point of $\mathcal{L}$ such that $\alpha\in
D_\lambda$. Then,
$f(\alpha)= G_\lambda(\alpha)P_\lambda(\alpha)e^{-Dp(\alpha)}=
P_\lambda(\alpha)=\omega(\alpha)$, as desired.

\section{Similar results in the disk}

The previous definitions and proofs can be adapted to produce analogous results
in the disk. To do so one just needs to replace the Euclidean distance used in
$\C$ by the pseudo-hyperbolic distance
\[
\rho(z,\zeta)=\left|\frac{z-\zeta}{1-\bar\zeta z}\right|\qquad z,\zeta\in\D ,
\]
and the Euclidean divided differences by their hyperbolic version
\[
\begin{split}
\delta^0 \omega(\lambda_1)  &=\omega(\lambda_1)\ ,\\
\delta^j\omega(\lambda_1,\ldots,\lambda_{j+1}) 
&=\displaystyle\frac{\Delta^{j-1}\omega(\lambda_2,\ldots,\lambda_{j+1})-\Delta^{
j-1}\omega(\lambda_1,\ldots,\lambda_j)}{\frac{\lambda_{j+1}-\lambda_1}{1-\bar
\lambda_1\lambda_{j+1}}}\qquad j\geq 1.\\
\end{split}
\]
In this context a function $\phi:\D\longrightarrow\R_+$ is a \emph{weight} if
\begin{itemize}
\item[(wd1)] There is a constant $K>0$ such that  $\phi(z)\ge K \ln\bigl(\frac
1{1-|z|}\bigr)$.
\item[(wd2)]  There are  constants  $D_0>0$ and $E_0>0$ such that whenever $\rho(z,\zeta)\le1 $ then
\[
\phi(z)\le D_0 \phi(\zeta)+E_0.
\]
\end{itemize}

The model for the associated spaces
\[
A_\phi=\{f\in H(\mathbb D) : \sup_{z\in\D}|f(z)|e^{-B\phi(z)}<\infty\ \textrm{for some $B>0$}\},
\]
is the Korenblum algebra $A^{-\infty}$, which corresponds to the choice
$e^{-\phi(z)}=1-|z|$. The interpolating sequences for this and similar algebras
have been characterised in \cite{Br-Pa} and \cite{Ma}.

With these elements, and replacing the factors $z-\alpha$ by
$\frac{z-\alpha}{1-\bar\alpha z}$ when necessary, we can follow the proofs above
and, mutatis mutandis, show that Theorem~\ref{main} also holds in this
situation. 

The only point that requires further justification is the validity of
Lemma~\ref{unif} for the weights $\phi$. Condition (i) is a standard consequence
of the open mapping theorem for (LF)-spaces applied to the restriction map
$\mathcal R_\Lambda$, and the same proof as in \cite[Lemma 2.2.6]{Br-Ga} holds.
Applying (i) to the sequences $\omega_\lambda(\Lambda)$ defined by
\[
\omega_\lambda(\lambda')=
\begin{cases}
1\quad&\textrm{if $\lambda'=\lambda$}\\
0\quad&\textrm{if $\lambda'\neq\lambda$}
\end{cases}
\]
we have functions $f_\lambda\in A_\phi$ interpolating these values and with
growth control independent of $\lambda$.
Since $1=|f_\lambda(\lambda)-f_\lambda(\lambda')|$, an estimate on the
derivative of $f_\lambda$ shows that
for some $C>0$ and $\varepsilon>0$ the pseudohyperbolic disks 
$D_H(\lambda,\varepsilon e^{-C\phi(\lambda)})=\{z\in\D
:\rho(z,\lambda)<e^{-C\phi(\lambda)}\}$ are pairwise disjoint. In particular the
sum of their areas is finite, hence
\[
\sum_{\lambda\in\Lambda}(1-|\lambda|)^2 e^{-2C\phi(\lambda)}<+\infty\ .
\]
From this and condition (wd1) we finally obtain (ii).

\end{document}